\documentclass[11pt]{article}

\usepackage{graphicx}
\usepackage{todonotes}

\def\beq{\begin{equation}}
\def\eeq{\end{equation}}
\def\beqa{\begin{eqnarray}}
\def\eeqa{\end{eqnarray}}
\def\beqs{\begin{eqnarray*}}
\def\eeqs{\end{eqnarray*}}

\newcommand{\rmd}{{\rm d}}

\def\EQ#1{Eq.~(\ref{#1})}

\def\dotprod{ \!\stackrel{\scriptscriptstyle\bullet}{{}}\! }

\newcommand{\UU}{{\mathcal U}}

\newcommand{\TT}{{\mathcal T}}
\newcommand{\BB}{{\mathcal B}}
\newcommand{\FF}{{\mathcal F}}

\newcommand{\HH}{{\mathcal H}}
\newcommand{\hh}{\tilde{{\mathcal H}}}
\newcommand{\utd}{\tilde{{\mathcal U}}}

\newcommand{\EE}{{\mathcal E}}
\newcommand{\GG}{{\mathcal G}}

\newcommand{\RR}{\mathbf{R}}
\newcommand{\uu}{\mathbf{u}}
\newcommand{\ta}{\mathbf{\tau}}
\newcommand{\n}{\mathbf{n}}

\newcommand{\pp}{\partial}

\begin{document}


\title{Boundary Integral Analysis for the Non-homogeneous \\ 3D Stokes Equation}

\author{
L.~J.~Gray\thanks{119 Berwick Drive, Oak Ridge, TN  37830 USA.},
Jas Jakowski\thanks{Oak Ridge High School, Oak Ridge, TN  37830}, 
M.~N.~J.~Moore\thanks{
Corresponding Author:
Department of Mathematics, Florida State University,
Tallahassee, Fl. \newline mnmoore2@fsu.edu},
and Wenjing Ye\thanks{
Department of Mechanical and Aerospace Engineering, The Hong Kong University of Science and Technology, Hong Kong} }
\maketitle

\abstract{A regular-grid volume-integration algorithm is developed for the non-homogeneous 3D Stokes equation.  Based upon the observation that the Stokeslet $\UU$ is the Laplacian of a function $\HH$,  the volume integral is reformulated as a simple boundary integral,  plus a remainder domain integral.  
The modified source term  in this remainder integral is everywhere zero on the boundary and can therefore be continuously extended as zero to a regular grid covering the domain. The volume integral can then be  evaluated on the grid.  
Applying this method to the Navier-Stokes equations will require obtaining velocity gradients, and thus an efficient algorithm for post-processing these derivatives is also discussed.  To validate the numerical implementation, test results employing a linear element Galerkin approximation are presented.}

\section{Introduction}
\label{Introduction}

Boundary integral analysis for fluids has primarily involved 
potential (inviscid irrotational)  and slow viscous (linear Stokes)
flow \cite{Poz}.   In particular, the break-up  
and coalescence  of free surface 
flows using integral formulations  of the
Laplace \cite{Eggers1,Leppinen,RT1,2Fluid,coal}
and Stokes \cite{Youngren,Rallison,SieroLister,Stone,PozJCP} equations
have a long and successful history.  
Other examples of Stokes flow applications include rigid particle suspensions \cite{af2014fast}, vesicle suspensions \cite{quaife2014high, liu2017dynamics, kabacaouglu2017low}, surfactant flows \cite{kropinski2011efficient}, electrohydrodynamics \cite{veerapaneni2016integral}, microswimmers \cite{kanevsky2010modeling, keaveny2011applying}, biological cell modeling \cite{Pozbio} including mitotic cell division \cite{nazockdast2017fast, wu2017forces}, viscous erosion \cite{mitchell2017generalized, Quaife2018}, and micro-electro-mechanical systems (MEMS)\cite{YeMEMS,Frangi}.

With a more general fluid equation, the conversion to  
a boundary integral statement cannot be completed exactly.  
To take an important example, the 
Navier-Stokes equation for a viscous incompressible fluid
\cite{Lady} on the domain $\Omega$ becomes the 
standard Stokes boundary integral together with
a volume term of the form 
\beq
\label{volterm}
 \int_{\Omega} \UU_{kj}(Q,P)\FF_j(Q)  \,{\rm d}\Omega_Q \,\, ,
\eeq
where $\UU$ is the Stokes Green's function and (with the summation convention
and $\rho_d$ the fluid density)
\beq
\label{NSF}
\FF_j =  \rho_d \uu_k \frac{\pp \uu_j}{\pp x_k} \,\, .
\eeq
is the non-linear convection term  \cite{Australia}.
Similarly, the
integral formulation for a non-Newtonian fluid must also contend with
a {\lq}pseudo body force{\rq} domain integral 
of the form in \EQ{volterm} \cite{nonnewton}.
For the moving boundary simulations mentioned above, 
constructing a finite element volume mesh
at every time step would be computationally demanding, if not impossible
considering the extreme geometries.
To overcome this difficulty, we present here a {\lq}body-fitted{\rq} volume discretization
for the non-homogeneous 3D Stokes equation.

Methods that have been employed for volume integral evaluation
include direct finite element mesh calculations
\cite{Australia,nonnewton,IngberV1,OfSteinbach}, and approximate
conversion to a boundary integral using either dual reciprocity
\cite{Power,ZhengColeman} or interior line integration \cite{Hematiyan,Fata_V}. 
The embedded boundary method in
\cite{Biros2004} and the fast Poisson solvers \cite{Askham,greengardV}
exploit an easily constructed regular grid covering the domain,
in conjunction with a Fast Multipole Method.

In  particular, the volume integration in  \cite{Askham} for the Poisson solver
and the work herein both rely on the continuous extension of a
function from $\Omega$ to the grid. 
In the former,  the source $\FF$ is extended by solving
an exterior Dirichlet boundary integral equation,
with boundary values given by $\FF$.  For the Stokes algorithm herein,
as well as the previous  Laplace
and elasticity treatments discussed in \cite{KYG,Jim,Dmitry,EFM2017},
a function $\FF^0$  is defined on $\Omega$ by means of
an {\it interior} boundary integral solution (again with boundary values from $\FF$).
With $\BB$ denoting the covering box, and with the help of Green's Theorem,
the modified body force integral becomes
\beq
\label{remainder}
 \int_{\BB} \UU_{kj}(Q,P)\left( \FF_j - \FF^0_j \right)(Q)  \,{\rm d}\Omega_Q \,\, .
\eeq
Here, $\FF-\FF^0$ is  defined as zero outside $\Omega$ and is therefore continuous on the entire grid.
The evaluation of this integral can be carried out by decomposing $\BB$ into (uniform) cuboid
cells, and using simple linear interpolation over {\it all} cells.  In other words,
with the continuous zero extension the
{\lq}partial cells{\rq} that straddle the boundary are treated as any
other cell and  the location of the domain boundary can effectively be ignored.
Moreover, as $ \FF_j - \FF^0_j \equiv 0$ on the boundary, the Green's function
singularity does not present any problem. 
Thus, unlike for the Poisson algorithm in \cite{Askham}, 
it is not necessary to refine the cells near
the boundary.
The calculations for $\FF^0$,  including interior point evaluations at
cell vertices, are however expensive; nevertheless, as in
\cite{Biros2004,Askham,greengardV,Tornberg2008},
fast methods can be employed \cite{Dmitry}.


A second key aspect of the work in
\cite{KYG,Jim,Dmitry,EFM2017}
is the construction (in a simple analytic form) of a function $\HH$ that satisfied
\beq
\EE(\HH) = \GG \,\, ,
\eeq
where $\EE$ is the partial differential equation and $\GG$ the corresponding
Green's function.   For the 3D Stokes equation, it will turn out that 
the function $\HH$ will satisfy the Stokes
equations with {\it zero pressure}.  That is,
\beqa
\label{HStokes}
\mu\nabla^2\HH  & = & \UU \\
\nabla\dotprod\HH & = & 0  \,\, , \nonumber
\eeqa 
where again $\UU$ is the Stokeslet, the Laplacian is applied componentwise
and the divergence is applied to the columns of $\HH$. 
The numerical treatment of the Stokes volume integral therefore
has much in common with the Laplace algorithm in \cite{KYG}, 
and this somewhat simplifies the numerical implementation.  

It is expected that this domain integral algorithm can form the basis
of an effective 
numerical solution of the two non-linear problems mentioned above, 
incompressible Navier-Stokes and non-Newtonian flows.  For the former,
the body force is given in \EQ{NSF}, and consequently
the computation of $\FF^0$ will require obtaining velocity
gradients on the boundary.   
The evaluation of these derivatives will therefore
be considered in Section \ref{VG}.
Although a non-linear solution algorithm will
not be considered herein, it is worth noting that previous papers
 \cite{Ding2,Yenonlinear} have successfully dealt with non-linear equations
in conjunction with the regular grid algorithm.

\section{$3$D Stokes' Equation}
\label{3D}

The equations for Stokes flow in three dimensions are \cite{Poz}
\beqa
\label{StokesF}
\mu\nabla^2\uu - \nabla p & = & \FF \nonumber \\
\nabla \dotprod \uu  &  = &   0   \,\, ,
\eeqa
where the problem domain is $\Omega$ and
$\{ \uu,p,\mu\}$ are fluid velocity, pressure, and
viscosity. The source function $\FF$ can be
a given body force and/or nonlinear terms stemming from, for example,
the Navier-Stokes or non-Newtonian equations.  

The Stokes boundary/volume integral statement is \cite{Poz,Australia,Biros2004}
\beqa
\label{ExtLimit}
& & \int_{\Sigma}  \bigg\lbrack \TT_{ijk}(Q,P)\uu_i(Q)\n_j(Q)
  -  \UU_{kj}(Q,P) \ta_j(Q) \bigg\rbrack
\,\rmd\Sigma_Q  =  \nonumber \hphantom{---------} \\
& &  \hphantom{---------------}
 \int_{\Omega} \UU_{kj}(Q,P)\FF_j(Q)  \,{\rm d}\Omega_Q \,\, ,
\eeqa
where $\Sigma$ is the boundary of
$\Omega$, $\n$ the {\it exterior} unit normal , $\ta(Q)$ the surface force
(normal component of the stress), $\UU_{kj}(Q,P)$
the Stokeslet (Green's function) and $\TT_{ijk}(Q,P)$ the corresponding
stresslet.  In this equation, it is assumed that the singular integrals are
defined as a limit as $P$ approaches the boundary from {\it outside}
the domain (this only comes into play for the $\TT$ kernel). For later use
(see Section \ref{VG}), there is also an interior limit equation,
(\cite{Poz}, Eq. (2.3.11))
\beqa
\label{IntLimit}
-\uu(P) & + & \int_{\Sigma}  \bigg\lbrack \TT_{ijk}(Q,P)\uu_i(Q)\n_j(Q)
  -  \UU_{kj}(Q,P) \ta_j(Q) \bigg\rbrack
\,\rmd\Sigma_Q  =  \nonumber \hphantom{-} \\
& &  \hphantom{--------}
 \int_{\Omega} \UU_{kj}(Q,P)\FF_j(Q)  \,{\rm d}\Omega_Q \,\, ,
\eeqa
the only difference being the explicit appearance of
the {\lq}free term{\rq} $\uu(P)$.

With 
$Q=\{q_{\ell}\}$, $P=\{p_{\ell}\}$,
$\RR=\{R_{\ell}\}=Q-P$ and $r=\Vert \RR \Vert$ the distance,
the fundamental solutions are given by \cite{Poz}
\beqa
\UU_{kj}(Q,P) & = & \frac{1}{8\pi\mu} \left[
\frac{\delta_{kj}}{r} + \frac{R_kR_j}{r^3}  \right]  \,\, , \nonumber \\
\TT_{ijk}(Q,P) & = & -\frac{3}{4\pi} \left[
\frac{R_iR_jR_k}{r^5}
 \right]
 \label{UandT}
\eeqa
where $\delta_{kj}$ is the Kronecker delta function.  

As with the Laplace and elasticity volume algorithms \cite{KYG,Dmitry}, we seek to
represent the kernel $\UU$ as a derivative in such a way that
Green's Theorem can be used to reformulate the volume integral in \EQ{ExtLimit}.
Based upon the
isotropic elasticity analysis in \cite{Dmitry}, and the fact that the Stokeslet is
the elasticity Green's function with Poisson ratio $\nu = 1/2$ and
shear modulus equal to viscosity, we define
the function $\HH$ as
\beqa
\label{HH}
\HH_{kj}(Q,P) & = & \frac{1}{8\pi\mu^2} \left[
{\delta_{kj}} r - \frac{1}{12} \frac{\pp^2 }{\pp q_k \pp q_j}r^3  \right]  \,\, ,
\nonumber \\
& = &  \frac{1}{32\pi\mu^2} \left[              
3r\delta_{kj} - \frac{R_kR_j}{r}  \right]   \,\, .
\eeqa
Direct calculation (see Appendix A) shows that
\beqa
\label{HU}
\mu\nabla^2 \HH_{kj} & = & \UU_{kj} \nonumber \\
\nabla \dotprod \HH & = & 0 \,\, ,
\eeqa
which says that $\HH$ satisfies the zero pressure non-homogeneous
Stokes equations, with the Stokeslet as the pseudo body force.
This is completely analogous to the  Green's functions representations 
for Laplace  and elasticity that were employed in \cite{KYG,Dmitry}.  
Fom this point the volume integral treatment 
follows as in these previous articles.   However, as
\EQ{HU} is componentwise,   the Stokes formulation more closely
resembles the scalar Laplace implementation than the vector elasticity. 

Briefly then, defining the vector $\FF^0$ by
\beq
\label{lapF0}
\nabla^2\FF^0 = 0
\eeq
with boundary conditions $\FF^0 \big\vert_\Sigma = 
\FF\big\vert_\Sigma$, and using \EQ{HU} we can write
the volume integral as
\beq
\label{Vol3}
\mu\int_{\Omega} \nabla^2\HH_{kj}(Q,P)\FF^0_j(Q)  \,{\rm d}\Omega_Q  + 
\int_{\Omega} \UU_{kj}(Q,P)\left(\FF_j - \FF^0_j \right) (Q)  
                 \,{\rm d}\Omega_Q 
\,\, .
\eeq
Invoking Green's Theorem and  \EQ{lapF0}, the first term becomes
the boundary integral
\beq
\label{BI} 
 \mu\int_{\Sigma} \left(
 \frac{\pp \HH_{kj}}{\pp \n}(Q,P) \FF^0_j(Q) - 
 \HH_{kj} (Q,P)\frac{\pp F^0_j}{\pp \n}(Q)
 \right)   \,{\rm d}\Sigma_Q  \,\, .
\eeq
The $\FF^0$ boundary flux is obtained from a boundary integral solution of \EQ{lapF0}.
Note that unlike their counterparts in \EQ{UandT},
$\HH(Q,P)$ and its normal derivative
\beq
\label{Hnorderiv}
 \frac{\pp \HH_{kj}}{\pp \n} =
  \frac{1}{32\pi\mu^2 r} \left(
    {3\n\dotprod R}  \delta_{kj} 
    - \left( {\n_kR_j + \n_jR_k} \right)
    +\n\dotprod R \frac{R_k R_j }{r^2}  \right)
\eeq
are rather innocous at $Q=P$. 
$\HH$  behaves as $r$, and is therefore
continuous, while the normal derivative is merely discontinuous.

\subsection{Numerics: Remainder Volume Integral}
\label{Numerics}

All boundary integrals in this work are approximated using a standard
Galerkin procedure with linear triangular elements. 
As just noted,   the kernel functions in \EQ{BI} do not
diverge at $Q=P$, and thus the treatment of \EQ{BI} is somewhat simplified.
In this case all integrals are 
computed numerically, except for the coincident integral involving the
normal derivative of $\HH$.  Techniques for the singular and non-singular
boundary integrals are well known, so this section will briefly discuss the 
remainder volume integral  \EQ{RemVol}.  Further details can be found
in previous papers \cite{KYG,Dmitry}.

For the volume term, the Galerkin form can be exploited by
interchanging of the order of integration  \cite{OfSteinbach}
\beqa
\label{RemVol}
\int_{\Sigma} \psi_k (P) \int_{\Omega} \UU_{kj}(Q,P)\left(\FF_j - \FF^0_j \right) (Q)  
                 \,{\rm d}\Omega_Q \, {\rm d}\Sigma_P  = \quad\quad   \quad\quad  \\
 \quad\quad  \quad\quad
 \int_{\Omega} \left(\FF_j - \FF^0_j \right) (Q)   \int_{\Sigma} \psi_k (P) \UU_{kj}(Q,P)
           \, {\rm d}\Sigma_P      \,{\rm d}\Omega_Q   \nonumber   \,\, .
\eeqa
The integrand in the boundary integral is known in analytic form 
and thus this integral can be handled partially analytically.  In the regular
grid approach, the points $Q$ will be the cell vertices, and the
interchange is especially
useful for the {\lq}near-singular{\rq} case when $Q$  is close to the boundary.
 
The other key aspect of \EQ{RemVol} is that by design
$\FF_j-\FF^0_j \equiv 0$ on $\Sigma$. This allows a continuous
extension of this function as zero outside the domain ,
and permits evaluation using a regular grid that covers the domain $\Omega$.
Using simple linear interpolation, the volume integration becomes simply a sum
over the vertex values.  The required interior 
values of $\FF^0$ can be computed from the 
interior point boundary integral equation once the
boundary $\FF^0$  flux has been obtained.  Again, the reader is asked
to consult \cite{KYG,Dmitry} for further details.

\subsection{$2$D Stokes' Equation}
\label{2D}

Although the
discussion herein is for three dimensions, for completeness we include
the expression for the 2D version of \EQ{HH}.  
The $\HH$ function  for 2D Stokes, denoted by $\hh$, 
is related to that for elasticity and,
as shown in \cite{Dmitry}, the general
form for  elasticity is (dropping the normalization constant)
\beq
\hh_{kj}  =
 \delta_{kj} r^2 \left(a + b\log(r) \right) + 
  \frac{\pp^2}{\pp q_k\pp q_j} \left\{ r^4 \left(c + d\log(r) \right) \right\} \,\, .
\eeq
In terms of the Poisson ratio, $\nu$,
the coefficients are 
 $b = \nu - 1$,
$d = (3 - 4\nu) / ( 64(1-\nu) )$ and $a$, $c$ are only required to satisfy 
$64c(1-\nu) + 2a(3-4\nu) = 8\nu^2-9\nu+2$.
For elasticity there is therefore leeway in the choice of these last two coefficients.

The Stokes equation corresponds to  $\nu=1/2$ and 
in this case $b=-1/2$, $d=1/32$, and
$16c+a=-1/4$.  However,  to satisfy  
$\nabla \dotprod \hh = 0$ and  $\mu\nabla^2\hh=\utd$,
where $\utd_{kj}$ is the two dimensional the Stokelet  
\cite{Poz} (Eq. $2.6.17$)
\beq
\utd_{kj}(Q,P)  =  
-\delta_{kj}\log(r) + \frac{R_kR_j}{r^2}  \,\, ,
 \label{2DU}
\eeq
it is necessary that  $c=-1/16$ and hence $a=3/4$.  Thus, in two dimensions
\beqa
\label{2DH}
\hh_{kj}  & = &  \frac{1}{4}  \delta_{kj}r^2\left\lbrack 3 - 2\log(r) \right\rbrack  +
  \frac{1}{32} \frac{\pp^2}{\pp q_k\pp q_j} r^4 \left\lbrack \log(r) -2 \right\rbrack
  \nonumber \\
  & = & \frac{1}{32} \bigg\lbrack \delta_{kj}r^2 \left( 17-12\log(r) \right) 
    + 2R_kR_j\left( 4\log(r) - 5\right) \bigg\rbrack  \,\, .
 \eeqa
A Maple code that  verifies this result can be found in  the Appendix.

\section{Test Calculations}
\label{Tests}

\subsection{$\FF=0$}

Before testing the evaluation of the volume integral, we first provide
evidence that the boundary integral in \EQ{ExtLimit} 
({\it i.e.}, without the volume integral) has been
implemented correctly. To this end, we solve the homogeneous 
Stokes equations interior and exterior to the unit sphere, with
velocity  boundary conditions on the upper half, $z\ge 0$,
and force boundary data on the lower half.  These boundary values,
as well as the corresponding exact solution,  are taken from
a point source placed exterior to the problem domain, specifically,
the first column of the Green's functions in \EQ{UandT}. 
Table \ref{P1Table}
lists the average square nodal  errors
\beq
\left[ \frac{1}{N} \sum_{k=1}^N \epsilon_k^2 \right]^{1/2}
\eeq
for each component of the computed velocity and force.  
Discretizations with $N_E=376$ and
$N_E=1504$ elements and two exterior points
$P_1=(-2.0,0,0)$ and $P_2=(1.5,0,0)$  were employed.  As expected, the
errors are smaller for the velocity, and there is 
reasonable decay in the errors for the finer mesh.  Moreover,  
as $P_2$ is closer to the
boundary, these errors are appropriately larger than for $P_1$.

\begin{table}[!ht] 
\label{P1Table} 
\begin{center}
\begin{tabular}[t]{|cc||lll||lll|}
\hline
\multicolumn{2}{|c||}{}  & \multicolumn{3}{c||}{$N_E=376$}  &
\multicolumn{3}{c|}{$N_E=1504$}  \\
\multicolumn{2}{|c||}{}  & \multicolumn{1}{c}{$x$}  &
\multicolumn{1}{c}{$y$}  & \multicolumn{1}{c||}{$z$}  &
 \multicolumn{1}{c}{$x$}  &
\multicolumn{1}{c}{$y$}  & \multicolumn{1}{c|}{$z$} \\
 \hline
 $P_1$  &  $\uu$  & 
1.094E-4  &   5.670E-5  &   3.634E-5  &
3.796E-5  &   1.409E-5  &   7.812E-6 \\
 &   $\ta$   & 
   5.520E-4  &   2.042E-4   &  2.347E-4  &
2.214E-4  &   7.840E-5  &   9.376E-5   \\
 \hline
$P_2$  &   $\uu$  &  
2.487E-4  &   8.597E-5  &   1.320E-4  &
6.050E-5  &  1.985E-5   &    2.009E-5  \\
 &   $\ta$   &   
1.515E-3   &  1.327E-3  &   1.136E-3  &
4.502E-4  &  2.751E-4  &   2.504E-4   \\
 \hline \hline	
\end{tabular}
\caption{ Mean square errors for boundary velocity and force for a homogenous
Stokes problem posed on the unit
sphere.  The mixed boundary conditions are from two exterior point source
locations,
$P_1=(-2.0,0,0)$ and $P_2=(1.5,0,0)$.}
\end{center}
\end{table}

Table \ref{P3Table} lists the corresponding errors for the exterior Stokes
problem, the two point sources located at $P_3=(0,0.7,0)$ 
and $P_4=(0,0,0.8)$.  The results again indicate that the Stokes boundary
integral equation has been implemented correctly.

\begin{table}[!ht] 
\label{P3Table} 
\begin{center}
\begin{tabular}[t]{|cc||lll||lll|}
\hline
\multicolumn{2}{|c||}{}  & \multicolumn{3}{c||}{$N_E=376$}  &
\multicolumn{3}{c|}{$N_E=1504$}  \\
\multicolumn{2}{|c||}{}  & \multicolumn{1}{c}{$x$}  &
\multicolumn{1}{c}{$y$}  & \multicolumn{1}{c||}{$z$}  &
 \multicolumn{1}{c}{$x$}  &
\multicolumn{1}{c}{$y$}  & \multicolumn{1}{c|}{$z$} \\
 \hline
$P_3$  &  $\uu$  & 
  1.964E-3  &   1.749E-3  &  3.412E-3  &
 4.472E-5  &   3.674E-5  &  4.184E-5  \\
 &   $\ta$   & 
  1.468E-2  &   1.530E-2  &   1.152E-2  &
  7.571E-4  &   6.693E-4  &   5.617E-4  \\
   \hline
$P_4$  &  $\uu$  &
  8.042E-4  &   7.977E-4  &   1.458E-3   &
  2.281E-5   &  1.991E-5   &  4.091E-5  \\
 &   $\ta$   &
 4.614E-2  &   1.222E-2  &   2.792E-2  &
 5.633E-3  &   2.772E-3 &    5.118E-3  \\
 \hline \hline	
\end{tabular}
\caption{ Mean square errors for boundary velocity and force for a homogenous
Stokes problem posed exterior to the unit
sphere.  The mixed boundary conditions are from two exterior point source
locations,
$P_3=(0,0.7,0)$ and $P_4=(0,0,0.8)$.}
\end{center}
\end{table}

\subsection{$\FF\neq 0$}

To verify the evaluation of the volume integral, \EQ{BI} and \EQ{RemVol}, 
a problem similar
to the above point source tests is employed.   Note that \EQ{HU},
\beqa
\label{HU2}
\mu\nabla^2 \HH_{kj} & = & \UU_{kj} \nonumber \\
\nabla \dotprod \HH & = & 0 \,\, ,
\eeqa
states that the columns of $\HH$ satisfy the Stokes equations, with
zero pressure and with forcing function $\FF$ given by the corresponding
column of $\UU$.   Thus, for an exterior point $P$, 
$\RR=Q-P$, $r=\Vert\RR\Vert$, and selected column $j$,
the velocity field is
\beq
       u_{kj} =    \frac{1}{32\pi\mu^2} r \left( 3\delta_{kj} - \frac{R_k}{r} \frac{R_j}{r}
       \right)    \,\, .
\eeq
The corresponding stress field for $\HH$ is
\beq
\sigma_{kjl}  = \frac{1}{16\pi\mu} \left( \frac{R_kR_lR_j}{r^3} +
  \frac{ R_l\delta_{kj} + R_k\delta_{jl} - 
 R_j\delta_{kl} }{r}  \right)
\eeq

\begin{table}[!ht] 
\label{HTable} 
\begin{center}
\begin{tabular}[t]{|cc||lll||lll|}
\hline
\multicolumn{2}{|c||}{}  & \multicolumn{3}{c||}{$N_E=376$}  &
\multicolumn{3}{c|}{$N_E=1504$}  \\
\multicolumn{2}{|c||}{}  & \multicolumn{1}{c}{$x$}  &
\multicolumn{1}{c}{$y$}  & \multicolumn{1}{c||}{$z$}  &
 \multicolumn{1}{c}{$x$}  &
\multicolumn{1}{c}{$y$}  & \multicolumn{1}{c|}{$z$} \\
 \hline
 $P_1$  &  $\uu$  & 
2.815E-4  &   3.495E-5  &   8.432E-5  &
7.382E-5  &   8.187E-6  &   1.908E-5 \\
 &   $\ta$   & 
6.332E-4  &   2.033E-4   &  3.076E-4  &
2.199E-4  &   7.934E-5  &   1.016E-4   \\
 \hline
$P_2$  &   $\uu$  &  
3.314E-4  &   3.644E-5  &   1.039E-4  &
8.111E-5  &  8.700E-6   &    2.392E-5  \\
 &   $\ta$   &   
6.721E-4   &  2.400E-4  &   3.604E-4  &
2.165E-4  &  8.181E-5  &   1.512E-4   \\
 \hline
$P_3$  &   $\uu$  &  
2.579E-4  &   2.300E-5  &    5.769E-4  &
6.686E-5  &  7.441E-6   &    1.336E-5  \\
 &   $\ta$   &   
5.881E-4   &  2.913E-4  &   2.790E-4  &
2.010E-4  &  9.461E-5  &    8.569E-5   \\
\hline
$P_4$  &   $\uu$  &  
1.643E-4  &   1.322E-5  &   5.353E-5  &
5.299E-5  &  3.347E-6   &   1.315E-5  \\
 &   $\ta$   &   
6.246E-4   &  1.984E-4  &   3.690E-4  &
2.043E-4  &   7.412E-5  &    1.100E-4  \\
 \hline \hline	
\end{tabular}
\caption{Mean square errors for boundary velocity and force for a nonhomogenous
Stokes problem posed on the unit
sphere.  The mixed boundary conditions are from the function $\HH(Q,P)$ 
with $P$ the exterior points 
$P_1=(2.0,0,0)$, $P_2=(1.5,0,0)$,
$P_3=(0,1.3,0)$, $P_4=(0,0,1.2)$.}
\end{center}
\end{table}

As with the point source test,
the domain is the unit sphere, and velocity is specified for
$z\ge 0$ and traction for $z<0$.  For the remainder volume integral,
the box that covers the sphere is $-1.1 < \{x,y,z\} < 1.1$, subdivided with
a $40\times40\times40$ grid.
Table \ref{HTable} lists the errors in the computed velocity
and traction for four exterior
points located at $P_1=(2.0,0,0)$ 
$P_2=(1.5,0,0)$,
$P_3=(0,1.3,0)$ and $P_4=(0,0,1.2)$.  
The results confirm that the volume terms, the $F_0$ boundary
integral equation and the remainder volume integral, 
have been implemented correctly.

\section{Velocity Gradients}
\label{VG}

An important potential application of the regular grid algorithm would be
a {\lq}boundary mesh only{\rq} solution of
the Navier-Stokes equations for viscous incompressible flow.
The pseudo body force $\FF$, \EQ{NSF}, is  a function of 
velocity gradients \cite{Australia,Power}, and thus the boundary
conditions for $\FF^0$ require knowledge of these derivatives on $\Sigma$
(gradient values will  also required at interior grid vertices, 
but herein we deal solely with the boundary).    In this section,  
the Stokes implementation of the gradient method 
proposed in \cite{GPK2004} will be discussed.

To start, it is a simple matter to differentiate the {\it interior} limit 
equation \EQ{IntLimit} to obtain an expression for the gradient components,
\beqa
\label{GradEQ}
- \uu_{k,m}(P) & = & \int_{\Sigma}  \bigg\lbrack \,
\UU_{kj,m}(Q,P) \ta_j(Q)   
 -  \TT_{kjl,m}(Q,P)\uu_i(Q)\n_j(Q) \, \bigg\rbrack \,\rmd\Sigma_Q 
\nonumber \\
&  + &   \int_{\Omega} \UU_{kj,m}(Q,P)\FF_j(Q)  \,{\rm d}\Omega_Q \,\, ,
\eeqa
where the subscript $,m$ indicates differentiation with respect to 
the coordinate $P_m$.   Note that once \EQ{ExtLimit} has been solved,
everything on the right hand side of \EQ{GradEQ} is known.
Nevertheless,  a simple implementation of this equation
is problematic: the computation is clearly very 
expensive and it involves the hypersingular kernel $ \TT_{kjl,m}(Q,P)$.
Moreover, given the approximations involved, 
\EQ{GradEQ} would probably not produce a highly accurate result.  
As a result, a
variety of alternative boundary integral gradient algorithms that have been considered,
\cite{ChatiMukh,GMEtangential,Guiggiani94,Schwab99,WildeFerri}
being just a partial list.
The papers \cite{Manticgrad3,ZhaoLan} and the references therein provide
a more complete overview of this topic. 

To address the computational cost, the  
key observation  in \cite{GPK2004} is that, unlike for the basic velocity 
and traction boundary integral equations,  the {\it exterior} limit gradient equation
\beqa
\label{GradEQext}
0 & = & \int_{\Sigma}  \bigg\lbrack \,
\UU_{kj,m}(Q,P) \ta_j(Q)   
 -  \TT_{kjl,m}(Q,P)\uu_i(Q)\n_j(Q) \, \bigg\rbrack \,\rmd\Sigma_Q 
\nonumber \\
&  + &   \int_{\Omega} \UU_{kj,m}(Q,P)\FF_j(Q)  \,{\rm d}\Omega_Q \,\, ,
\eeqa
and the interior limit  are not the same equation.   
Thus, subtracting the two equations, interior minus exterior, yields a new and useful
expression for the gradient.   The only terms that do not cancel in the limit difference
are those that are discontinuous crossing
the boundary, and this  immediately eliminates all non-singular boundary integrals. Moreover,  {\it in the volume integral}, the Stokeslet singularity is sufficiently weak 
that this term is also continuous across the boundary, and we can therefore write
\beq
\label{GradEQLIM}
\uu_{k,m}(P) = Lim_{\leftrightarrow}  \int_{\Sigma}  \bigg\lbrack \,
\TT_{kjl,m}(Q,P)\uu_i(Q)\n_j(Q) -
\UU_{kj,m}(Q,P) \ta_j(Q)  \, \bigg\rbrack
\,\rmd\Sigma_Q  \,\, ,
\eeq
where $Lim_{\leftrightarrow} $ indicates the limit difference. 

A Galerkin implementation is employed to deal with the 
second issue, the hypersingular kernel $ \TT_{ijk,m}(Q,P)$. 
With Galerkin, \EQ{GradEQLIM}
reduces  to the coincident and adjacent edge singular integrations for
the $\TT$ integral, while for the 
weaker $\UU_{,m}$ singularity  it is solely the coincident
integral that contributes.   Although the integration work is minimal,
there is  an additional computational expense:
note that the Galerkin implementation of \EQ{GradEQLIM}
couples the nodal values of a gradient component. As a consequence,
for each of the
nine  components  $\uu_{k,m}$, $1\le k,m \le 3$ on the boundary,
the solution of an $N\times N$ linear system is required,
$N$ the number of boundary nodes.  However, the
coefficient matrix is the same for each component, and moreover
it is sparse, symmetric positive definite.  Thus, for moderate sized
problems, only one matrix factorization is required, while for large
scale problems an efficient iterative solver can be employed.

Differentiating \EQ{UandT} with respect to the coordinates of $P$ the
derivative kernels are 
\beqa
\UU_{kj,m}(Q,P) & = & \frac{1}{8\pi\mu} \left[
\delta_{kj} \frac{R_m}{r^3}  - \delta_{km}\frac{R_j}{r^3}
                 -\delta_{jm}\frac{R_k}{r^3} + 3\frac{R_kR_jR_m}{r^5}
\right]  \,\, , \nonumber \\
\TT_{kjl,m}(Q,P) & = & -\frac{3}{4\pi} \left[
- \delta_{km} \frac{R_jR_l}{r^5}  
- \delta_{jm} \frac{R_kR_l}{r^5}
- \delta_{lm} \frac{R_kR_j}{r^5}  \right. \nonumber \\
&+&   \left.  5\frac{R_kR_jR_lR_m}{r^7}  \right]  \,\, .
 \label{hhderiv}
\eeqa
The coincident and adjacent edge singular integrals are evaluated
as in \cite{GPK2004,3dgalerkin}, and  the details will not
be repeated here.  However, it is worth noting that much of the
calculation is exact: the boundary limit
is handled analytically, and the singular integrals are
computed partially analytically.   Importantly, 
this allows the exact cancellation of
a potentially divergent term (of the form $1/\epsilon$, $\epsilon\rightarrow 0$ 
the distance to the boundary) that arises in the coincident $\TT_{,m}$ integration.
Further details about the singular integration can be found
in the cited references.

\subsection{Gradient Tests}

To verify the implementation of the  algorithm, complete velocity
gradients $\uu_{k,m}$ have computed for several  of the test 
problems discussed in Section \ref{Tests}.  
Table \ref{GradTable} lists the mean square 
gradient errors for (a) the interior homogeneous problem with point source 
at $P=(2,0,0)$;  (b)  the exterior homogeneous problem with point source 
at $P=(0,0.7,0)$; and (c) the non-homogeneous interior problem with
$P=(0,0,1.2)$.   For the last example 
the exact solution is the $\HH$ derivatives
\beq
\HH_{kj,m}(Q,P)  =  \frac{1}{32\pi\mu^2} \left[
-3\delta_{kj} \frac{R_m}{r}  + \delta_{km}\frac{R_j}{r}
                 +\delta_{jm}\frac{R_k}{r} - \frac{R_kR_jR_m}{r^3}
\right]  \,\, .
\eeq
with $j=1$ .  The results 
demonstrate that the volume integral does in fact cancel out of the limit
difference equation.

\begin{table}[!ht] 
\label{GradTable} 
\begin{center}
\begin{tabular}[t]{|c||lll||lll|}
\hline
\multicolumn{1}{|c||}{}  & \multicolumn{3}{c||}{$N_E=376$}  &
\multicolumn{3}{c|}{$N_E=1504$}  \\
 \hline
\multicolumn{1}{|c||}{$m$}  & 
\multicolumn{1}{c}{$\uu_{1,m}$}  & \multicolumn{1}{c}{$\uu_{2,m}$} & 
\multicolumn{1}{c||}{$\uu_{3,m}$}  &
\multicolumn{1}{c}{$\uu_{1,m}$}  & \multicolumn{1}{c}{$\uu_{2,m}$} & 
\multicolumn{1}{c|}{$\uu_{3,m}$}   \\
 \hline
 $1$  &  
1.002E-4    &     6.311E-5    &     6.833E-5   &
1.817E-5     &     1.199E-5    &     1.400E-5   \\
 $2$  & 
5.547E-5    &    5.718E-5    &     3.893E-5   &
9.240E-6    &    1.116E-5    &     7.820E-6    \\
  $3$  &   
6.335E-5    &     4.120E-5    &     5.443E-5  &   
9.940E-6    &     8.053E-6    &     1.168E-5 \\
\hline
  $1$ &  
4.428E-4    &     6.620E-4    &     3.751E-4  &
3.957E-5    &     5.462E-5    &     3.926E-5   \\
$2$ & 
 1.046E-3    &     5.296E-4    &     4.392E-4  &
 4.657E-5    &     3.758E-5    &     3.472E-5   \\
  $3$  &  
 5.526E-4    &     4.901E-4    &     4.443E-4  &
 2.818E-5    &     3.620E-5    &     2.649E-5   \\
 \hline
   $1$ &  
 1.135E-4    &     6.056E-5    &     6.881E-5   &
 2.187E-5    &     1.251E-5     &     1.435E-5  \\
 $2$ &   
 4.324E-5    &     5.383E-5    &     3.878E-5  &
 9.990E-6    &     1.094E-5    &     8.033E-6   \\
  $3$   &     
 5.934E-5    &     3.758E-5    &     5.362E-5  &
 1.198E-5     &     7.676E-6    &     1.114E-5   \\
 \hline \hline	
\end{tabular}
\caption{ Average square errors for boundary gradient components for 
(a)  interior homogeneous problem with source $P=(2,0,0)$;
(b) exterior  homogenous with source $P=(0,0.7,0)$;
and (c) non-homogeneous interior with $P=(0,0,1.2)$ for the unit
sphere.}
\end{center}
\end{table}

At first it might appear strange that the gradient errors are smaller
than for the corresponding computed velocity solution.  However, recall 
that these tests are mixed boundary value problems, with velocity boundary
conditions specified 
on half the sphere.  Naturally, the gradient calculation is strongly dependent 
on the input surface velocity, and thus in these examples the algorithm
is  working with exact data on half the boundary.

\section{Conclusion}
\label{Conclusion}

A regular grid volume integration algorithm
for the non-homogeneous 3D Stokes equation has been presented.
The key to modifying the original volume integral
is to represent the Green's function (Stokeslet) as the Laplacian
of  a function $\HH$.  This is analogous to previous volume integral treatments
for the Laplace  \cite{KYG,Jim} and elasticity equations  \cite{Dmitry,EFM2017}, 
as  the Laplacian can be viewed the zero-pressure Stokes equations.  
With the function $\HH$, the domain integral exactly transforms to a simple
boundary integral, plus a  volume term wherein the modified source function
is everywhere zero on the boundary.  The continuous zero extension of this source
allows this volume integral to be computed on a regular grid of cells covering
the domain.   

An effective Stokes volume integral technique allows the possibility
of treating nonlinear equations, {\it e.g.},  Navier-Stokes and 
non-Newtonian equations.  Nonlinear analyses with the regular grid approach
have been previously carried out \cite{Ding2,Yenonlinear}, and as well
there have been nonlinear integral equation solutions obtained with 
other volume methods \cite{Biros2004,Australia,Power,nonnewton}.
It is therefore reasonable to expect that this work will lead to an efficient
integral equation algorithm for  the nonlinear fluids noted above. 
Regarding the Navier-Stokes solver, the grid algorithm requires
boundary values of the surface gradient, and  it has been shown herein
that this post-processing calculation can be executed efficiently.

It would be of interest to repeat the successful potential flow studies
of coalescence \cite{coal} and Rayleigh-Taylor break-up \cite{RT1,2Fluid}
with the more complex fluids mentioned above.
The Laplace calculations employ
cylindrical coordinates $\{\rho,\theta,z\}$, and with $\theta$ integrated out
they become two-dimensional $\{\rho,z\}$ analyses.
The axi-symmetry also requires that any body force is independent 
of the polar angle, and in this situation it is straightforward to implement the
regular grid method  \cite{KYG}. For Stokes however,
axi-symmetry means that the body force vector, rather than being 
independent of $\theta$,  rotates
properly with $\theta$ \cite{Poz}.  As an example, the surface normal 
for an axi-symmetric geometry is 
not independent of $\theta$, but rather of the form
$\n=\left(n_x\cos(\theta),n_y\sin(\theta),n_z\right)$.
With $\FF$ of this form, the 
volume algorithm does not immediately carry over, and a modified
algorithm must be developed.  We hope to return to this issue
in the future.

As a second application, we hope to combine the framework developed here with the ideas laid out in \cite{weakcoupling} to simulate general viscoelastic flows. In this context, a polymeric stress field develops in conjunction with the flow and serves as the body force in the non-homogeneous Stokes system. Given an initial polymeric stress, the boundary-integral framework could be used to determine the corresponding flow, which would then be used to update the  stress field
at the next time step. Such a method would thus capture the nonlinear 
feedback between the polymeric stress and the flow with accuracy and efficiency.

\bibliography{Stokes}

\section*{Acknowledgements}

The authors are grateful to Prof. B. Quaife for kindly pointing out 
several important references.
L. J. Gray gratefully acknowledges an MTS Visiting Professorship Grant 
at the University of Minnesota, and he
would like to thank Profs. S. Mogilevskaya, J. Labuz and S. Crouch for the hospitality
at the Department of Civil Engineering.
The participation of  J.  Jakowski was facilitated by J. Williams and
Dr. D. Pickel of the Mathematics
Department, Oak Ridge High School. 
M. N. J. Moore acknowledges support from the Simons Foundation Collaboration Grants for Mathematicians, award ID 524259.

\section*{Appendix A}
\label{verifyH}

The Maple codes below can be used to  confirm \EQ{HU} for 3D and 2D.
That is, they establish that the
Stokes equations, with zero pressure, applied to $\HH_{kj}$ and $\hh_{kj}$
yield the corresponding 
Stokes Green's function, and these functions satisfy conservation of mass.

\subsection*{$\HH$ for $3D$}

\begin{verbatim}

q  := array(1..3);   p := array(1..3);   r  := array(1..3);
lapH  := array(1..3);  divH  := array(1..3);
EQ  := array(1..3,1..3);  H  := array(1..3,1..3);
G   := array(1..3,1..3);  del := array(1..3,1..3):
##
for k from 1 to 3 do
 for j from 1 to 3 do
  del[k,j] := 0;
 od;
  del[k,k] := 1;
od;
##
r[1] := q[1] - p[1];
r[2] := q[2] - p[2];
r[3] := q[3] - p[3];
rr   := sqrt( r[1]*r[1] + r[2]*r[2] + r[3]*r[3] );
##
##   Stokeslet
##
for k from 1 to 3 do
 for j from 1 to 3 do
   G[k,j]  :=  del[k,j]/rr + r[k]*r[j]/rr^3;
   G[k,j]  := G[k,j] / (8*Pi*mu);
 od;
od;
##
##  Elasticity function
##
a := 2 - 2*nu;    b := -(3-4*nu)/(24*(1-nu));    nu := 1/2;
for k from 1 to 3 do
 for j from 1 to 3 do
   H[k,j]  := a*rr*del[k,j] + b*diff(rr^3,q[k],q[j]); 
   H[k,j]  := H[k,j] / (8*Pi*mu^2);
 od;
od;
##
##   loop over columns of H
##
for j from 1 to 3 do

  lapH[1] := diff(H[1,j],q[1],q[1]) + diff(H[1,j],q[2],q[2]) + 
             diff(H[1,j],q[3],q[3]);
  lapH[2] := diff(H[2,j],q[1],q[1]) + diff(H[2,j],q[2],q[2]) + 
             diff(H[2,j],q[3],q[3]);
  lapH[3] := diff(H[3,j],q[1],q[1]) + diff(H[3,j],q[2],q[2]) + 
             diff(H[3,j],q[3],q[3]);

  lapH[1] := factor(normal(lapH[1]));
  lapH[2] := factor(normal(lapH[2]));
  lapH[3] := factor(normal(lapH[3]));

  divH[j] := diff(H[1,j],q[1]) +  diff(H[2,j],q[2]) + diff(H[3,j],q[3]);
  divH[j] := factor(normal(divH[j]));
  
##
##  Stokes' equation  mu*Lap(u) = grad(p) with p = 0
##  Stokes(H) = G
##

EQ[1,j] := factor(normal(expand( mu*lapH[1] - G[1,j] )));
EQ[2,j] := factor(normal(expand( mu*lapH[2] - G[2,j] ))); 
EQ[3,j] := factor(normal(expand( mu*lapH[3] - G[3,j] )));

od;    #   end column loop


\end{verbatim}

\subsection*{$\hh$ for $2D$}

\begin{verbatim}

R := array(1..2);

del[1,1] := 1;
del[2,2] := 1;
del[1,2] := 0;
del[2,1] := 0;

R[1]:=x[1]-p[1]:
R[2]:=x[2]-p[2]:
rsq := R[1]*R[1] + R[2]*R[2]:
r := sqrt(rsq):

nu := 1/2;
b := nu - 1;
d := (3 - 4*nu) / ( 64*(1-nu) );
a := -1/4 -16*c;
c := -1/16;

for k from 1 to 2 do
 for j from 1 to 2 do
  H[k,j]:=  rsq*( a + b*ln(rsq)/2 ) * del[k,j] +
            diff( rsq^2*( c + d*ln(rsq)/2 ),x[k],x[j]):
   LapH[k,j] := diff(H[k,j],x[1],x[1]) + diff(H[k,j],x[2],x[2]):
   LapH[k,j] := subs( (2*x[1]-2*p[1])^2 =  4*(x[1]-p[1])^2,
                      (2*x[1]-2*p[1])^4 = 16*(x[1]-p[1])^4,
                      (2*x[2]-2*p[2])^2 =  4*(x[2]-p[2])^2,
                      (2*x[2]-2*p[2])^4 = 16*(x[2]-p[2])^4,LapH[k,j]);
   G[k,j]    :=  -del[k,j]*ln(rsq)/2 + R[k]*R[j]/rsq;
   chk[k,j]  := normal(expand( LapH[k,j] - G[k,j] ));
 od;
    DivH[k] := diff(H[k,1],x[1]) + diff(H[k,2],x[2]):
    DivH[k] := normal(expand( DivH[k] ));
od;
\end{verbatim}

\end{document}